\documentclass[reqno,10pt]{amsart}
\usepackage{amscd,amssymb,verbatim}

\theoremstyle{plain}

\newcommand\alp{\alpha}         

\newcommand\gam{\gamma}         
         
\newcommand\eps{\varepsilon}
\newcommand\zet{\zeta}
\newcommand\tet{\theta}

         \newcommand\Sig{\Sigma}

\newcommand\ome{\omega}         \newcommand\Ome{\Omega}

\newcommand\calB{{\mathcal{B}}}

\newcommand\calI{{\mathcal{I}}}

\newcommand\RR{\mathbb{R}}

\newcommand\ZZ{\mathbb{Z}}

\newcommand\CC{\mathbb{C}}

 \newcommand\gro{{\mathfrak{o}}}

\newcommand\nek{,\ldots,}
\newcommand\sdp{\times \hskip -0.3em {\raise 0.3ex
\hbox{$\scriptscriptstyle |$}}} 

\newcommand\Det{\operatorname{Det}}

\newcommand\IM{\operatorname{Im}}

\newcommand\MOD{\operatorname{mod}}

\newcommand\Ker{\operatorname{Ker}}

\newcommand\rank{\operatorname{rank}}
\newcommand\rk{\operatorname{rk}}

\newcommand\RE{\operatorname{Re}}

\renewcommand{\>}{\rangle}
\newcommand{\<}{\langle}

\theoremstyle{plain}
\newtheorem{Thm}{Theorem}
\newtheorem{Cor}{Corollary}
\newtheorem{Lem}{Lemma}
\newtheorem{Prop}{Proposition}
\newtheorem{Conjec}{Conjecture}

\newtheorem{Def}{Definition}

\theoremstyle{remark}

\newtheorem{Rem}{Remark}

\def\TeXref#1{%
        \leavevmode\vadjust{\setbox0=\hbox{{\tt
                \  {\tiny \textrm #1}}}%
        \theight=\ht0
        \advance\theight by \lineskip
        \kern -\theight \vbox to
        \theight{\rightline{\rlap{\box0}}%
        \vss}%
        }}%
\newif\ifShowLabels
\ShowLabelstrue
\newdimen\theight
\def\TeXrefEq#1{%
        \leavevmode\vadjust{\setbox0=\hbox{{\tt
                \  {\tiny \textrm #1}}}%
        \theight=\ht1
        \advance\theight by \lineskip
        \kern -\theight \vbox to
        \theight{\rightline{\rlap{\box0}}%
        \vss}%
        }}%

\newcommand{\refss}[1]{\S\ref{SS:#1}}

\newcommand{\refd}[1]{Definition ~\ref{D:#1}}

\newcommand{\refe}[1]{\eqref{E:#1}}

\newenvironment{thm}[1]%
        { \begin{Thm} \label{T:#1}  \ifShowLabels \TeXref{T:#1} \fi }%
        { \end{Thm} }

\renewcommand{\th}[1]{\begin{thm}{#1}  }
\renewcommand{\eth}{\end{thm} }

\newenvironment{lemma}[1]%
        { \begin{Lem} \label{L:#1}  \ifShowLabels \TeXref{L:#1} \fi }%
        { \end{Lem} }

\newcommand{\lem}[1]{\begin{lemma}{#1} }
\newcommand{\elem}{\end{lemma}}

\newenvironment{propos}[1]%
        { \begin{Prop} \label{P:#1}  \ifShowLabels \TeXref{P:#1} \fi }%
        { \end{Prop} }

\newcommand{\prop}[1]{\begin{propos}{#1} }
\newcommand{\eprop}{\end{propos}}

\newenvironment{corol}[1]%
        { \begin{Cor} \label{C:#1}  \ifShowLabels \TeXref{C:#1} \fi }%
        { \end{Cor} }
\newcommand{\cor}[1]{\begin{corol}{#1}  }
\newcommand{\ecor}{\end{corol}}

\newenvironment{conjec}[1]%
        { \begin{Conjec} \label{Conj:#1}  \ifShowLabels \TeXref{C:#1} \fi }%
        { \end{Conjec} }
\newcommand{\conj}[1]{\begin{conjec}{#1}  }
\newcommand{\econj}{\end{conjec}}

\newenvironment{defeni}[1]%
        { \begin{Def} \label{D:#1}  \ifShowLabels \TeXref{D:#1} \fi }%
        { \end{Def} }
\newcommand{\defe}[1]{\begin{defeni}{#1}  }
\newcommand{\edefe}{\end{defeni}}

\newenvironment{remark}[1]%
        { \begin{Rem} \label{R:#1}  \ifShowLabels \TeXref{R:#1} \fi }%
        { \end{Rem} }
\newcommand{\rem}[1]{\begin{remark}{#1}}
\newcommand{\erem}{\end{remark}}

\newcommand{\eq}[1]%
        { \ifShowLabels \TeXrefEq{E:#1} \fi
           \begin{equation} \label{E:#1} }
\newcommand{\eeq}{\end{equation}}
\newcommand{\meq}[1]%
        { \ifShowLabels \TeXrefEq{E:#1} \fi
           \begin{multline} \label{E:#1} }
\newcommand{\emeq}{\end{multline}}

\newcommand{\prf}{ \begin{proof} }
\newcommand{\eprf}{ \end{proof} }
\newcommand{\Label}[1]{\label{#1}  \ifShowLabels \TeXref{#1} \fi }

\ShowLabelsfalse

\newcommand{\n}{\nabla}
\newcommand{\na}{\nabla_\alp}

\renewcommand{\b}{\bullet}

\newcommand{\even}{\operatorname{even}}

\newcommand{\Arg}{\operatorname{\mathbf{Arg}}}

\newcommand{\Mon}{\operatorname{Mon}}

\newcommand{\TRS}{T^{\operatorname{RS}}}
\newcommand{\TTur}{T^{\operatorname{comb}}}

\newcommand{\Detgrtet}{\Det_{{\operatorname{gr},\tet}}}

\newcommand{\p}{\pi_1(M)}
\newcommand{\Rep}{\operatorname{Rep}(\p,\CC^n)}
\newcommand{\Repo}{\operatorname{Rep_0}(\p,\CC^n)}
\newcommand{\Reph}{\operatorname{Rep}^u(\p,\CC^n)}
\newcommand{\Repho}{\operatorname{Rep}^u_0(\p,\CC^n)}

\newcommand{\ccomp}{C}
\begin{document}
\begin{flushright}
\end{flushright}

\title[]{A Refinement of the Ray-Singer  Torsion}
\author[]{Maxim Braverman$^\dag$}
\address{Department of Mathematics\\
        Northeastern University   \\
        Boston, MA 02115 \\
        USA
         }
\email{maximbraverman@neu.edu}
\author[]{Thomas Kappeler$^\ddag$}
\address{Institut fur Mathematik\\
         Universitat Z\"urich\\
         Winterthurerstrasse 190\\
         CH-8057 Z\"urich\\
         Switzerland
         }
\email{tk@math.unizh.ch}
\thanks{${}^\dag$Supported in part by the NSF grant DMS-0204421.\\
\indent${}^\ddag$Supported in part by the Swiss National Science Foundation.}
\begin{abstract}
We propose a refinement of the Ray-Singer torsion, which can be viewed as an analytic counterpart of the refined combinatorial torsion
introduced by Turaev. Given a closed, oriented manifold of odd dimension with fundamental group $\Gamma$, the refined torsion is a complex
valued, holomorphic function defined for representations of $\Gamma$ which are close to the space of unitary representations.  When the
representation is unitary the absolute value of the refined torsion is equal to the Ray-Singer torsion, while its phase is determined by the
$\eta-$invariant. As an application we extend and improve a result of Farber about the relationship between the absolute torsion of
Farber-Turaev and the $\eta$-invariant.

\smallskip \noindent {\bf R\'esum\'e}.\ {\bf Raffinement de la torsion de Ray-Singer.}\  Nous proposons un raffinement de la torsion analytique de Ray-Singer,
qui peut etre consider\'e comme un \'equivalent analytique du raffinement de la torsion combinatoire introduit par Turaev. Soit $M$ une
vari\'et\'e ferm\'ee et orient\'ee de dimension impaire et de groupe fondamental $\Gamma$. La torsion analytique raffin\'ee est une fonction
holomorphe \`a valeurs complexes, d\'efinie pour les repr\'esentations de $\Gamma$, qui sont proches de l'espace des repr\'esentations
unitaires. Dans le cas o\`u la repr\'esentation est unitaire, la valeur absolue de la torsion analytique raffin\'ee est \'egale \`a la torsion
de Ray-Singer d\`es que sa phase est d\'etermin\'ee par invariant $\eta$. Comme application, nous  g\'en\'eralisons et am\'eliorons un resultat
de Farber concernant la relation entre la torsion absolue de Farber-Turaev et invariant $\eta$.
\end{abstract}
\maketitle

\subsection*{Version fran\c caise abr\'eg\'ee}\Label{S:introd}
Soit $E\to M$ un fibr\'e complexe plat o\`u $M$ est une vari\'et\'e ferm\'ee et orient\'ee de dimension impaire. Pour un ensemble ouvert de
connections $\n$ plates et acycliques de $E$, qui  contient toutes les connections hermitiennes et acycliques, nous proposons  un raffinement de
la torsion analytique  de Ray-Singer $\TRS(\n)$. D\`es que $\TRS(\n)$ est un nombre positive r\'eel, $T= T(\n)$ est dans le cas g\'en\'eral un
nombre complexe et, comme cons\'equence, a une phase non-triviale. Ce raffinement peut \^etre consid\'er\'e comme un \'equivalent analytique du
raffinement du concept de la torsion combinatoire introduit par Turaev \cite{Turaev86,Turaev90} et d\'evelopp\'e plus tard par Farber et Turaev
\cite{FarberTuraev99,FarberTuraev00}. Bien que, en g\'en\'eral,  $T$ ne soit pas \'egal \`a la torsion de Turaev, les deux torsions sont
\'etroitement li\'ees.

Si $\dim{M}\equiv 1\ (\MOD 4)$ ou si le rang du fibr\'e $E$ est divisible par 4, la torsion analytique raffin\'ee $T=T(\n)$ est ind\'ependente de
tous les choix faits pour la d\'efinir. Si $\dim{M}\equiv 3\ (\MOD 4)$, $T(\n)$ d\'epend du choix  d'une vari\'et\'e compacte et orient\'ee $N$, le
bord orient\'e de celle-ci est diff\'eomorphe \`a deux copies disjointes de $M$, mais seulement \`a un  facteur $i^{k\cdot\rk{E}}$ \ ($k\in \ZZ$)
pr\`es.

Si la connection $\n$ est hermitienne, c'est \`a dire, s' il existe une m\'etrique hermitienne pour $E$ qui est invariante par $\n$, alors la
torsion analytique raffin\'ee $T$ est un nombre complexe dont la valeur absolue est \'egale \`a la torsion de Ray-Singer et dont la phase est
d\'etermin\'ee par invariant $\eta$ de l'op\'erateur de signature. Si $\n$ n'est pas hermitienne, les relations entre la torsion analytique
rafffin\'ee, la torsion de Ray-Singer et invariant $\eta$ sont un peu plus compliqu\'ees, cf. \refss{IlocalRS}.

Une des  propri\'et\'es  les plus importantes de la torsion analytique raffin\'ee est qu'elle d\'epend, dans un sens appropri\'e , d'une
mani\`ere {\em holomorphe} de la connection $\n$. Le fait qu'il soit possible d'utiliser la torsion de Ray-Singer et invariant $\eta$ pour
d\'efinir
 une fonction holomorphe permet d'appliquer les m\'ethodes de l'analyse complexe pour l'\'etude des deux invariants. En
particulier, nous obtenons  une relation entre la torsion analytique raffin\'ee et la torsion combinatoire raffin\'ee de Turaev qui
g\'en\'eralise le th\'eor\`eme c\'el\`ebre de Cheeger-M\"uller concernant l'\'egalit\'e entre la torsion analytique de Ray-Singer et la torsion
combinatoire \cite{Cheeger79,Muller78}. Comme application, nous g\'en\'eralisons et am\'eliorons un resultat de Farber relatif \`a la relation
entre la torsion absolue de Farber-Turaev et invariant $\eta$.

Notre construction de la torsion analytique raffin\'ee utilise le superd\'eterminant r\'egularis\'e de  $\zeta$ de l'op\'erateur de signature.
On remarquera que cet op\'erateur est non autoadjoint.

L'id\'ee de repr\'esenter la torsion de Ray-Singer comme valeur absolue d'une function holomorphe d\'efinie pour l'espace des connections est
due \`a Burghelea et Haller, \cite{BurgheleaHaller_Euler}. En r\'eponse \`a une version pr\'eliminaire de notre article, Dan Burghelea nous a
signal\'e son projet (men\'e avec Haller), \cite{BurgheleaHaller_function}, de construction d'une fonction holomorphe utilisant des determinants
regularis\'es de certain op\'erateurs  non autoadjoint du type d'un op\'erateur Laplacien.

\subsection{The odd signature operator}\Label{SS:Ioddsign}
Let $M$ be a closed oriented manifold of odd dimension $\dim{}M=n= 2r-1$ and let $E$ be a complex vector bundle over $M$ endowed with a flat
connection $\n$. Let $\Ome^\b(M,E)$ denote the space of smooth differential forms on $M$ with values in $E$ and set
\[
   \Ome^{\even}(M,E)\ = \ \bigoplus_{p=0}^{r-1}\Ome^{2p}(M,E),
\]
Fix a Riemannian metric $g^M$ on $M$ and let $*:\Ome^\b(M,E)\to \Ome^{n-\b}(M,E)$ denote the Hodge $*$-operator. The (even part of) the {\em odd
signature operator} is the operator
 \(
     B_{\even}=B_{\even}(\n,g^M):\,\Ome^{\even}(M,E)\ \to\ \Ome^{\even}(M,E),
 \)
whose value on a form $\ome\in \Ome^{2p}(M,E)$ is defined by the formula
\begin{equation} \Label{E:Ioddsign}\notag
    B_{\even}\,\ome \ := \ i^r(-1)^{p+1}\,\big(\, *\n-\n*\,\big)\,\ome \ \in \ \Ome^{n-2p-1}(M,E)\oplus\Ome^{n-2p+1}(M,E).
\end{equation}
The odd signature operator was introduced by Atiyah, Patodi, and Singer, \cite[p.~44]{APS1}, \cite[p.~405]{APS2}, and, in the more general
setting used here, by Gilkey, \cite[p.~64--65]{Gilkey84}.

The operator $B_{\even}$ is an elliptic differential operator, whose leading symbol is symmetric with respect to any Hermitian metric $h^E$ on
$E$.

In this paper we define the refined analytic torsion in the case when the pair $(\n,g^M)$ satisfies the following simplifying assumptions. The
general case will be addressed elsewhere.

\subsection*{Assumption~I} The connection $\n$ is acyclic, i.e.,
\[
    \IM\big(\n|_{\Ome^{k-1}(M,E)}\big) \ =\  \Ker\big(\n|_{\Ome^k(M,E)}\big),\qquad \text{for every}\quad k=0\nek n.
\]
\subsection*{Assumption~II} The odd signature operator $B_{\even}=B_{\even}(\n,g^M)$ is bijective.

\bigskip
Note that all acyclic Hermitian connections satisfy Assumptions~I and II. By a simple continuity argument, these two assumptions are then satisfied
for all flat connections in an open neighborhood (in $C^0$-topology) of the set of acyclic Hermitian connections.

\subsection{Graded determinant}\Label{SS:Igrdet}
Set
\eq{Igrading}
  \Ome^k_+(M,E) \ := \ \Ker\,(\n\,*)\,\cap\,\Ome^k(M,E), \qquad
  \Ome^k_-(M,E) \ := \ \Ker\,(*\,\n)\,\cap\,\Ome^k(M,E).
\end{equation}
Assumption~II implies that $\Ome^k(M,E)\ = \ \Ome^k_+(M,E)\oplus \Ome^k_-(M,E)$. Hence, \refe{Igrading} defines a {\em grading} on $\Ome^k(M,E)$.

Define $\Ome^{\even}_\pm(M,E)= \bigoplus_{p=0}^{r-1}\Ome^{2p}_{\pm}(M,E)$ and let $B_{\even}^\pm$ denote the restriction of $B_{\even}$ to
$\Ome^{\even}_\pm(M,E)$. It is easy to see that $B_{\even}$ leaves the subspaces $\Ome^{\even}_\pm(M,E)$ invariant and it follows from
Assumption~II that the operators $B^\pm_{\even}:\Ome^{\even}_\pm(M,E)\to \Ome^{\even}_\pm(M,E)$ are bijective.

One of the central objects of this paper is the {\em graded determinant} of the operator $B_{\even}$. To construct it we need to choose a {\em
spectral cut} along a ray $R_{\tet}= \big\{\rho{}e^{i\tet}:0\le \rho<\infty\big\}$, where $\tet\in [-\pi,\pi)$ is an Agmon angle for $B_{\even}$.
Since the leading symbol of $B_{\even}$ is symmetric, $B_{\even}$ admits an Agmon angle $\tet\in (-\pi,0)$. Given such an angle $\tet$, observe that
it is an Agmon angle for $B_{\even}^\pm$ as well. The graded determinant of $B_{\even}$ is the non-zero complex number defined by the formula
\eq{Igrdeterminant}
   \Detgrtet(B_{\even}) \ := \ \frac{\Det_\tet(B_{\even}^+)}{\Det_\tet(B_{\even}^-)}.
\end{equation}
By standard arguments, $\Detgrtet(B_{\even})$ is independent of the choice of the Agmon angle $\tet\in (-\pi,0)$.

\subsection{A convenient choice of the Agmon angle}\Label{SS:IAgmon}
For $\calI\subset \RR$ we denote by $L_{\calI}$ the solid angle
\eq{ILab}\notag
    L_{\calI} \ = \  \big\{\, \rho e^{i\tet}:\, 0 < \rho<\infty,\, \tet\in \calI\,  \big\}.
\end{equation}

Though many of our results are valid for any Agmon angle $\tet\in (-\pi,0)$, some of them are easier formulated if the following conditions are
satisfied:

\begin{enumerate}
\item[\textbf{(AG1)}] \ $\tet\in (-\pi/2,0)$, and
\item[\textbf{(AG2)}] \ there are no eigenvalues of the operator $B_{\even}$ in the solid angles $L_{(-\pi/2,\tet]}$ and $L_{(\pi/2,\tet+\pi]}$.
\end{enumerate}

For the sake of simplicity of exposition, we will assume that $\tet$ is chosen so that these conditions are satisfied. Since the leading symbol of
$B_{\even}$ is symmetric, such a choice of $\tet$ is always possible.

\subsection{Relationship with the Ray-Singer torsion and the $\eta$-invariant}\Label{SS:Ietainv}
For a pair $(\n,g^M)$ satisfying Assumptions~I and II set
\eq{IcalB}
    \xi \ = \  \xi(\n,g^M,\tet) \ := \ \frac12\,\sum_{k=0}^{n-1}\,(-1)^k\,\zet_{2\tet}' \big(\,0,\,
    (-1)^{k+1}\,{(*\,\n)^2}{\big|_{\Ome^k_+(M,E)}}\,\big),
\end{equation}
where $\zet_{2\tet}'\big(\,s,\,(-1)^{k+1}\,{(*\,\n)^2}{\big|_{\Ome^k_+(M,E)}}\,\big)$ is the derivative with respect to $s$ of the $\zet$-function of
the operator ${(*\,\n)^2}{\big|_{\Ome^k_+(M,E)}}$ corresponding to the spectral cut along the ray $R_{2\tet}$, and $\tet$ is an Agmon angle
satisfying (AG1)-(AG2).

Let $\eta=\eta(\n,g^M)$ denote the $\eta$-invariant of the operator $B_{\even}(\n,g^M)$, cf. Definition~4.2 of \cite{BrKappelerRAT}. Note that,
since the operator $B(\nabla,g^{TM})$ is not, in general, self-adjoint, $\eta$ might be a complex number, cf. \cite{Gilkey84} and Section~4 of
\cite{BrKappelerRAT}. Theorem~7.2 of \cite{BrKappelerRAT} implies that,
\eq{IDetB-eta}
    \Detgrtet(B_{\even}) \ = \ e^{\xi(\n,g^M,\tet)}\cdot e^{-i\pi\eta(\n,g^M)}.
\end{equation}
This representation of the graded determinant turns out to be very useful, e.g., in computing the metric anomaly of $\Detgrtet(B_{\even})$.

If the connection $\n$ is Hermitian, then \refe{IcalB} coincides with the well known expression for the logarithm of the Ray-Singer torsion
$\TRS= \TRS(\n)$. Hence, for a Hermitian connection $\n$ we have
\[
    \xi(\n,g^M,\tet) \ = \ \log\,\TRS(\n).
\]
If $\n$ is not Hermitian but is sufficiently close (in $C^0$-topology) to an acyclic Hermitian connection, then Theorem~8.2 of \cite{BrKappelerRAT}
states that
\eq{IDetB-TRS}
    \log \TRS(\n) \ = \ \RE\,\xi(\n,g^M,\tet).
\end{equation}

Combining \refe{IDetB-TRS} and \refe{IDetB-eta}, we get
\eq{IDetB-TRSherm}
    \big|\,  \Detgrtet(B_{\even}) \,\big| \ = \ \TRS(\n)\cdot e^{\pi\IM\eta(\n,g^M)}.
\end{equation}
If $\n$ is Hermitian, then the operator $B_{\even}$ is self-adjoint and $\eta= \eta(\n,g^M)$ is real. Hence, for the case of an acyclic Hermitian
connection we obtain from \refe{IDetB-TRSherm}
\eq{IDetB-TRSherm2}\notag
    \big|\,  \Detgrtet(B_{\even}) \,\big| \ = \ \TRS(\n).
\end{equation}

\subsection{Definition of the refined analytic torsion}\Label{SS:Irefinedtorsion}
The graded determinant of the odd signature operator is not a differential invariant of the connection $\n$ since, in general, it depends on the
choice of the Riemannian metric $g^M$. Using the metric anomaly of the $\eta$-invariant computed by Gilkey, \cite{Gilkey84}, we investigate in
\S9 of \cite{BrKappelerRAT} the {\em metric anomaly} of the graded determinant and then use it to ``correct" the graded determinant and
construct a differential invariant -- the refined analytic torsion.

Suppose an acyclic connection $\n$ is given. We call a Riemannian metric $g^M$ on $M$ {\em admissible for $\n$} if the operator $B_{\even}=
B_{\even}(\n,g^M)$ satisfies Assumption~II of \refss{Ioddsign}. We denote the set of admissible metrics by $\calB(\n)$. The set $\calB(\n)$ might be
empty. However, admissible metrics exist for all flat connections in an open neighborhood (in $C^0$-topology) of the set of acyclic Hermitian
connections, cf. Proposition~6.8 of \cite{BrKappelerRAT}.

\defe{reftorsion}
The {\em refined analytic torsion} $T(\n)$ corresponding to an acyclic connection $\n$, satisfying $\calB(\n)\not=\emptyset$, is defined as follows:
fix an admissible Riemannian metric $g^M\in \calB(\n)$ and let  $\tet\in (-\pi,0)$ be an Agmon angle for $B_{\even}(\n,g^M)$.
\begin{enumerate}
\item
If $\dim M\equiv 1\ (\MOD 4)$ then
\[
   T(\n)\ =\ T(M,E,\n)\ :=\ \Detgrtet\big(B_{\even}(\n,g^M)\big) \ \in \ \CC\backslash{0}.
\]
\item
If $\dim M\equiv 3 (\MOD 4)$ choose a smooth compact oriented manifold $N$ whose oriented boundary is diffeomorphic to two disjoint copies of $M$
(since $\dim{}M$ is odd, such a manifold always exists, cf. \cite{Wall60}). Then
\begin{equation}\Label{E:Irefantor3}\notag
    T(\n) \ = \ T(M,E,\n,N) \ := \
        \Detgrtet(B_{\even})\cdot \exp\Big(\,i\pi\,\frac{\rank E}2\,  \int_{N}\, L(p)\,\Big) \ \in \ \CC\backslash{0},
\end{equation}
where $L(p)$ is the Hirzebruch $L$-polynomial in the Pontrjagin forms of a Riemannian metric on $N$ which is a product near $M$.
\end{enumerate}
\edefe
Note that $\int_{N}\, L(p)$ is real and, hence, $|T(\n)| = \ |\Detgrtet(B_{\even})|$.

If\/ $\n$ is close enough to an acyclic Hermitian connection, then $\calB(\n)\not=\emptyset$ and it is shown in \S\S11-12 of \cite{BrKappelerRAT},
that {\em $T(\n)$ is independent of the choices of the admissible metric $g^M$ and the Agmon angle $\tet\in (-\pi,0)$}. However, if $\dim M\equiv 3
(\MOD 4)$, then the refined analytic torsion {\em  does depend on the choice of the manifold $N$}. The quotient of the refined torsions corresponding
to different choices of $N$ is a complex number of the form $i^{k\cdot\rank{E}}$ $(k\in \ZZ)$. Hence, if $\rank{}E$ is even then $T(\n)$ is well
defined up to a sign, and if $\rank{}E$ is divisible by 4, then $T(\n)$ is a well defined complex number.


\subsection{Comparison with the Ray-Singer torsion}\Label{SS:IlocalRS}
The equality \refe{IDetB-TRSherm} implies that, if $\n$ is $C^0$-close to an acyclic Hermitian connection, then
\eq{ITRS-T}
    \log\,\frac{|T(\n)|}{\TRS(\n)} \ = \ \pi\,\IM\,\eta(\n,g^M)
\end{equation}
In particular, if $\n$ is an acyclic Hermitian connection, then
 \(
    \big|\,T(\n)\,\big| \ = \ \TRS(\n).
 \)

Let $\Arg_\n$ denote the unique cohomology class $\Arg_\n\in H^1(M,\CC/\ZZ)$ such that for every closed curve $\gam\in M$ we have \/
 \(
    \det\big(\,\Mon_\n(\gam)\,\big) = \exp\big(\, 2\pi i\<\Arg_\n,[\gam]\>\,\big),
 \)
where $\Mon_\n(\gam)$ denotes the monodromy of the flat connection $\n$ along the curve $\gam$ and $\<\cdot,\cdot\>$ denotes the natural pairing
 \(
    H^1(M,\CC/\ZZ)\,\times\, H_1(M,\ZZ) \to  \CC/\ZZ.
 \)
Then, cf. Theorem~12.8 of \cite{BrKappelerRAT}, if $\n$ is $C^0$-close to an acyclic Hermitian connection, then
\eq{IRaySinger3}
    \log\,\frac{|T(\n)|}{\TRS(\n)} \ \ = \ \ \pi\,\big\<\, [L(p)]\cup \IM\Arg_{\n},[M]\,\big\>.
\end{equation}

If \/ $\dim M\equiv 3\ (\MOD\ 4)$, then $L(p)$ has no component of degree $\dim{M}-1$ and, hence,
 \(
    |T(\n)|  =  \TRS(\n).
  \)

\subsection{The refined analytic torsion as a holomorphic function on the space of representations}\Label{SS:Iholomorphic}
One of the main properties of the refined analytic torsion $T(\n)$ is that, in an appropriate sense, it depends holomorphically on the connection.
Note, however, that the space of connections is infinite dimensional and one needs to choose an appropriate notion of a holomorphic function on such
a space. As an alternative one can view the refined analytic torsion as a holomorphic function on a finite dimensional space, which we shall now
explain.

The set $\Rep$ of all $n$-dimensional complex representations of $\p$ has a natural structure of a complex algebraic variety. Each
representation $\alp\in \Rep$ gives rise to a vector bundle $E_\alp$ with a flat connection $\na$, whose monodromy is isomorphic to $\alp$. Let
$\Repo\subset \Rep$ denote the set of all representations $\alp\in \Rep$ such that the connection $\na$ is acyclic. Further we denote by
$\Reph\subset \Rep$ the set of all unitary representations and set $\Repho= \Reph\cap\Repo$.

Denote by $V\subset \Repo$ the set of representations $\alp$ for which there exists a metric $g^M$ so that the odd signature operator
$B_{\even}(\n,g^M)$ is bijective (i.e., Assumption~II of \refss{Ioddsign} is satisfied). It is easy to see that $V$ is an open neighborhood of
the set $\Repho$ of acyclic unitary representations.

For every $\alp\in V$ one defines the refined analytic torsion $T_\alp:= T(\na)$. Corollary~13.11 of \cite{BrKappelerRAT} states that the function
$\alp\mapsto T_\alp$ is holomorphic on the open set of all non-singular points of $V$.

\subsection{Comparison with Turaev's torsion}\Label{SS:ITuraev}
In \cite{Turaev86,Turaev90}, Turaev introduced a refinement $\TTur_\alp(\eps,\gro)$ of the combinatorial torsion associated to a representation
$\alp$ of $\pi_1(M)$. This refinement depends on an additional combinatorial data, denoted by $\eps$ and called the {\em Euler structure} as well as
on the {\em cohomological orientation} of $M$, i.e., on the orientation $\gro$ of the determinant line of the cohomology $H^\b(M,\RR)$ of $M$. There
are two versions of the Turaev torsion -- the homological and the cohomological one. It is more convenient for us to use the cohomological Turaev
torsion as it is defined in Section~9.2 of \cite{FarberTuraev00}. For $\alp\in \Repo$ the cohomological Turaev torsion $\TTur_\alp(\eps,\gro)$ is a
non-vanishing complex number.

Theorem~10.2 of \cite{FarberTuraev00} computes the quotient of the Turaev and the Ray-Singer torsions. Combined with \refe{IRaySinger3} this leads to
the following  result: for every Euler structure $\eps$ and every cohomological orientation $\gro$, there exists an open neighborhood $V'\subset V$
of $\Repho$ such that for every $\alp\in V'$
\eq{IRelogTTTur}
   \left|\,\frac{T_\alp}{\TTur_\alp(\eps,\gro)}\,\right| \ = \ \left|\, f_{\eps,\gro}(\alp)\,\right|,
\end{equation}
where $f_{\eps,\gro}(\alp)$ is a holomorphic function of $\alp\in V'$, which is explicitly calculated in \S14.4 of \cite{BrKappelerRAT}.

Let $\Sig$ denote the set of singular points of the complex analytic set $\Rep$. The refined analytic torsion $T_\alp$ is a non-vanishing holomorphic
function of $\alp\in V\backslash\Sig$. By the very construction \cite{Turaev86,Turaev90,FarberTuraev00} the Turaev torsion is a non-vanishing
holomorphic function of $\alp\in \Repo$. Hence, $(T_\alp/{\TTur_\alp})^2$ is a holomorphic function on $V'\backslash\Sig$. If the absolute values of
two non-vanishing holomorphic functions are equal on a connected open set then the functions must be equal up to a factor $\phi\in \CC$ with
$|\phi|=1$. This observation and \refe{IRelogTTTur} lead to the following generalization of the Cheeger-M\"uller theorem, cf. Theorem~14.5 of
\cite{BrKappelerRAT}: \ {\em For each connected component $\ccomp$ of \/ $V'$, there exists a constant $\phi_\ccomp= \phi_\ccomp(\eps,\gro)\in \RR$,
depending on $\eps$ and $\gro$, such that
\eq{IlogTTTur}
    \frac{T_\alp}{\TTur_\alp(\eps,\gro)} \ = \ e^{i\phi_\ccomp}\,f_{\eps,\gro}(\alp).
\end{equation}}

It would be interesting to extend \refe{IlogTTTur} to a generalization of the Bismut-Zhang extension \cite{BisZh92} of the Cheeger-M\"uller
theorem.

\subsection{Application: Phase of the Turaev torsion}\Label{SS:Farber}
The equality \refe{IlogTTTur} provides, in particular, a relationships between the phases of the Turaev torsion and the refined analytic
torsion. If $\alp$ is a unitary representation, then \refe{IDetB-eta} and \refd{reftorsion} express the phase of the refined analytic torsion in
terms of the $\eta$-invariant. Thus we obtain a relationship between the phase of the Turaev torsion and the $\eta$-invariant for the case of a
unitary representation $\alp$, see Theorem~14.9 of \cite{BrKappelerRAT} for the precise form of this relationship. In particular, this leads to
an extension and a new proof of a theorem of Farber \cite{Farber00AT} about the relationship between the sign of the absolute torsion of
Farber-Turaev, \cite{FarberTuraev99}, and the $\eta$-invariant, cf. Theorem~14.12 of  \cite{FarberTuraev99}.

\providecommand{\bysame}{\leavevmode\hbox to3em{\hrulefill}\thinspace} \providecommand{\MR}{\relax\ifhmode\unskip\space\fi MR }
\providecommand{\MRhref}[2]{%
  \href{http://www.ams.org/mathscinet-getitem?mr=#1}{#2}
} \providecommand{\href}[2]{#2}

\end{document}